\documentclass[12pt,twoside]{amsart}
\usepackage{amssymb}

\nonstopmode

\numberwithin{equation}{section} 

\font\tengothic=eufm10 scaled\magstep 1
\font\sevengothic=eufm7 scaled\magstep 1
\newfam\gothicfam
          \textfont\gothicfam=\tengothic
          \scriptfont\gothicfam=\sevengothic

\newcommand{\s}{\; | \;}
\newcommand{\fall}{\mbox{for all} ~}
\newcommand{\mif}{\mbox{if} ~}

\newcommand{\Z}{\mathbb{Z}}
\newcommand{\R}{\mathbb{R}}

\DeclareMathOperator{\chara}{char}

 \DeclareMathOperator{\Tor}{Tor}

\newcommand{\fm}{{\mathfrak m}}

\DeclareMathOperator{\pnt}{\raise 0.5mm \hbox{\large\bf.}}

\newcommand \pp[1] {^{\langle #1 \rangle}}
\newcommand \pl[1] {_{[ #1 ]}}
\newcommand \dd[1] {_{\langle #1 \rangle}}

\newcommand{\WLP}{Weak Lefschetz property}

\newtheorem{theorem}{Theorem}[section]
\newtheorem{lemma}[theorem]{Lemma}
\newtheorem{proposition}[theorem]{Proposition}
\newtheorem{corollary}[theorem]{Corollary}
\newtheorem{conjecture}[theorem]{Conjecture}

\newtheorem{prop-def}[theorem]{Proposition and Definition}

\theoremstyle{definition}
\newtheorem{definition}[theorem]{Definition} 
\newtheorem{remark}[theorem]{Remark}
\newtheorem{example}[theorem]{Example}

\newtheorem{notation}[theorem]{Notation}

\begin{document}

\title[Empty simplices and graded Betti numbers]
{Empty simplices of polytopes and graded Betti numbers}

\author[Uwe Nagel]{Uwe Nagel}

\address{Department of Mathematics,
University of Kentucky, 715 Patterson Office Tower,
Lexington, KY 40506-0027, USA}

\email{uwenagel@ms.uky.edu}


\begin{abstract}
The conjecture of Kalai, Kleinschmidt, and Lee on the number of empty simplices of a
simplicial polytope is established by relating it to the first graded Betti numbers of
the polytope. The proof allows us to derive explicit optimal bounds on the number of
empty simplices of any given dimension. As a key result, we prove optimal bounds for the
graded Betti numbers of any standard graded $K$-algebra in terms of its Hilbert function.
\end{abstract}

\maketitle

\section{Introduction}

Let $P \subset \R^d$ be a simplicial $d$-polytope, i.e.\ the $d$-dimensional convex hull
of finitely many points in $\R^d$ such that all its faces are simplices. The simplest
combinatorial invariant of $P$ is its $f$-vector $\underline{f} = (f_{-1}, f_0,\ldots,
f_{d-1})$ where $f_{-1} := 1$ and $f_i$ is the number of $i$-dimensional faces of $P$ if
$i \geq 0$. McMullen conjectured in \cite{McM-70} a characterization of the possible
$f$-vectors. In order to state his conjecture we use an equivalent set of invariants, the
$h$-vector $\underline{h} := (h_0,\ldots, h_s)$. It is defined as the sequence of
coefficients of the polynomial
\[
\sum_{j=0}^s h_j z^j := \sum_{j=0}^d f_{j-1}\cdot z^j (1-z)^{d-j}.
\]
The $f$-vector can be recovered from the $h$-vector because
$$
f_{j-1} = \sum_{i = 0}^j \binom{d-i}{j-i} h_i.
$$
Using $h$-vectors we can state McMullen's conjecture which has become a proven statement
by combining the results of  Billera and Lee \cite{BL2} and  Stanley \cite{stanley2}
(cf.\ also \cite{McM-simple}).

\begin{theorem}[$g$-theorem] \label{g-thm}
A sequence $\underline{h} = (h_0,\dots,h_s)$
of positive integers is the $h$-vector of a simplicial $d$-polytope if and only
if $s = d$ and $\underline{h}$ is an SI-sequence, i.e.\
$\underline{h}$ satisfies:
\begin{itemize}
\item[(i)] (Dehn-Sommerville equations) $h_i = h_{d-i}$ for $i = 0,\ldots,d$;
\item[(ii)] $\underline{g} := (h_0, h_1 - h_0,\ldots,h_{\lfloor
  \frac{d}{2} \rfloor} - h_{\lfloor \frac{d}{2} \rfloor - 1})$ is an
  O-sequence.
\end{itemize}
\end{theorem}

Being an O-sequence is a purely numerical condition (cf.\ Section \ref{sec-Betti}). Note
that O-sequences are precisely the Hilbert functions of Artinian standard graded
$K$-algebras.

In order to prove sufficiency of these conditions, in \cite{BL2} Billera and Lee
construct, for each SI-sequence $\underline{h} := (h_0,\ldots, h_d)$, a certain
simplicial $d$-polytope $P_{BL}(\underline{h})$ whose  $h$-vector is the given
SI-sequence $\underline{h}$. The Billera-Lee polytopes are rather particular which has
lead to expectations that they have some extremal properties. In order to state one such
instance recall (cf.\ \cite{Kalai-02}) that an {\em
  empty simplex}  of the polytope $P$ is a smallest subset $S$ of the
vertex set of $P$ such that $S$ is not a face of $P$, but each proper subset of $S$ is a
face of $P$. Sometimes, empty simplices are called {\em
  missing faces}. They are just minimal
non-faces of the vertex set of $P$. Empty simplices play an important role in the
classification of polytopes (cf., e.g, \cite{Kalai-87} and Remark \ref{rem-skeleta}). In
\cite{Kalai-94}, Kalai states as Conjecture 2:

\begin{conjecture}[Kalai, Kleinschmidt, Lee] \label{conj}
For all simplicial $d$-polytopes with prescribed $h$-vector $\underline{h}$, the number
of $j$-dimensional empty simplices is maximized by the Billera-Lee polytope
$P_{BL}(\underline{h})$.
\end{conjecture}

Kalai has pointed out in \cite{Kalai-02}, Theorem 19.5.35, that this conjecture is a
consequence of results in \cite{MN3}, but his argument needs some adjustment. The
starting point of this note is to give a detailed proof of this conjecture which is
established in Theorem \ref{thm-conj}.

The construction of the Billera-Lee polytopes is rather involved. In general, the number
of empty $j$-simplices of a given Billera-Lee polytope $P_{BL} (\underline{h})$ has not
been known. Hence, the proof of Conjecture \ref{conj} leaves open the problem of giving
an explicit bound in terms of the $h$-vector. The bulk of this paper is devoted to
solving this problem. The key is given by our proof of Conjecture \ref{conj}. It
identifies the number of missing $j$-simplices of the polytope $P$ with a certain graded
Betti number of its Stanley-Reisner ring $K[P]$. Since the $h$-vector of $P$ is
determined by the Hilbert function of $K[P]$, we are lead to consider the problem of
finding sharp upper bounds for the graded Betti numbers of the Stanley-Reisner ring
$K[P]$ in terms of its Hilbert function. We solve this problem in Section \ref{sec-Betti}
in greater generality, namely for Gorenstein algebras with the Weak Lefschetz property
(Theorem \ref{thm-gor-betti}). Its proof requires   explicit bounds for all graded Betti
numbers of any standard graded $K$-algebra $A$ in terms of its Hilbert function. These
are established in Theorem \ref{thm-Betti-bounds}. They are optimal. Because of the
importance of graded Betti numbers, it seems fair to expect that Theorem
\ref{thm-Betti-bounds} will find applications in other contexts as well.

In Section \ref{sec-missing}, we apply the results of Section \ref{sec-Betti} to derive
explicit optimal bounds for the number of missing $j$-simplices of a simplicial polytope
in terms of its $g$-vector (cf.\ Corollary \ref{cor-generators}). Note that the
$g$-vector is easily obtained from the $h$-vector (Definition \ref{deg-g-vector}). We
conclude with some applications. In particular, we bound the number of empty faces of
dimension $\leq k$ of a simplicial $d$-polytope in terms of $k$ and $f_0 - d$ (Corollary
\ref{cor-kalai-2-7}). Following Kalai \cite{Kalai-94}, such a bound is the key to a
central result of Perles \cite{perles} in the theory of arbitrary polytopes with ``few
vertices'' (cf.\ Remark \ref{rem-skeleta}). Finally, we show that very little information
on the $g$-vector is sufficient to bound the number of empty $j$-simplices of a
simplicial $d$-polytope if $d$ is large enough (Corollary \ref{cor-j-simp}). This result
slightly corrects and improves \cite{Kalai-94}, Theorem 3.8.

\section{The conjecture of Kalai, Kleinschmidt, and Lee} \label{sec-conj}

The goal of this section is to prove Conjecture \ref{conj}. To this end we need some more
notation. Let $P$ be a simplicial $d$-polytope. Denote its vertex set
 by $\{v_1,\ldots,v_{f_0}\}$ and let $R := K[x_1,\ldots,x_{f_0}]$ be the
polynomial ring in $f_0$ variables over an arbitrary field $K$. Then the {\em
Stanley-Reisner ring} of $P$ is $K[P] := R/I_P$ where the {\em Stanley-Reisner ideal} is
generated by all square-free monomials $x_{i_1} x_{i_2} \cdots x_{i_t}$ such that
$\{v_{i_1}, v_{i_2},\ldots,v_{i_t}\}$ is not a face of $P$. It is well-known (cf.\
\cite{BH-book}, Corollary 5.6.5) that $K[P]$ is a Gorenstein ring of dimension $d = \dim
P$. Since $h_1 = f_0 - d$, its minimal graded free resolution is of the form
\[
0 \rightarrow \bigoplus_{j \in \Z} R(-j)^{\beta_{h_1, j}^K (P)}
\rightarrow \dots
\rightarrow
\bigoplus_{j \in \Z} R(-j)^{\beta_{1, j}^K (P)} \rightarrow R
\rightarrow R/I \rightarrow
0.
\]
The non-negative integers $\beta_{i, j}^K (P) = \dim_K [\Tor_i^R (K[P],
  K)]_j$, $i, j \in \Z$, are called the {\em graded Betti numbers} of
$P$.

The following result is  shown in \cite{MN3}:

\begin{theorem} \label{thm-betti}
Let $K$ be a field of characteristic zero and let $P$ be a simplicial
$d$-polytope with $h$-vector $\underline{h}$.  Then we have for all
integers $i, j$:
$$
\beta_{i, j}^K (P) \leq \beta_{i, j}^K (P_{BL} (\underline{h})).
$$
\end{theorem}

\begin{proof}
The claim is a consequence of \cite{MN3}, Theorem 9.6, because its
proof shows (cf.\ page 57) that the extremal polytope that is not
specified in part
(b) of this theorem is indeed the Billera-Lee polytope $P_{BL}
(\underline{h})$.
\end{proof}

\begin{remark} The assumption on the characteristic of the field $K$ is
  needed to ensure that the Stanley-Reisner ring $K[P]$ has the
  so-called Weak Lefschetz property (cf.\ Section \ref{sec-Betti}). This property
  also plays a crucial role in Stanley's necessity part of the
  $g$-theorem in \cite{stanley2}.
\end{remark}

The Conjecture of Kalai, Kleinschmidt, and Lee follows now easily.

\begin{theorem} \label{thm-conj}
For all simplicial polytopes with prescribed $h$-vector $\underline{h}$, the number of
$j$-dimensional empty simplices is maximized by the Billera-Lee polytope
$P_{BL}(\underline{h})$.
\end{theorem}

\begin{proof}
It follows from its definition that ${\beta_{1, j}^K} (P)$ is the
number of minimal generators of degree $j$ of the Stanley-Reisner
ideal $I_P$. Since a $j$-dimensional empty face of $P$ corresponds to
a minimal generator of $I_P$ with degree $j+1$,  the Conjecture of
Kalai, Kleinschmidt, and Lee is a consequence of Theorem
\ref{thm-betti} applied with $i=1$.
\end{proof}

The combinatorial interpretation of the first Betti numbers allows us to drop the
assumption on the characteristic in Theorem \ref{thm-betti} for certain Betti numbers.

\begin{corollary} \label{first-betti}
Let $P$ be a simplicial $d$-polytope with $h$-vector $\underline{h}$.  Then we have for
all integers $j$:
$$
\beta_{1, j}^K (P) \leq \beta_{1, j}^K (P_{BL} (\underline{h})),
\quad \beta_{h_1 - 1, j}^K (P) \leq \beta_{h_1 - 1, j}^K (P_{BL}
(\underline{h})),
$$
and
$$
\beta_{h_1, j}^K (P) = \left \{ \begin{array}{ll}
0 & \mif j \neq h_1 + d \\
1 & \mif j = h_1 + d
\end{array} \right.
$$
\end{corollary}

\begin{proof}
Denote by $n_j (P)$ the number of empty $j$-simplices of $P$. We have seen that, for
every field $K$:
$$
n_{j-1} (P) = \beta_{1, j}^K (P).
$$
Let $K$ be a field of characteristic zero. Then Theorem \ref{thm-conj} provides
$$
n_{j-1} (P) \leq n_{j-1} (P_{BL} (\underline{h})).
$$
Let now $K$ be an arbitrary field. Then, applying the above equality again, the claim for
the first Betti numbers follows.

Since $K[P]$ is a Gorenstein ring, its minimal free resolution is self-dual. In
particular, for all integers $i, j$, we have
$$
\beta_{i, j}^K (P) = \beta_{h_1 - i,  h_1 + d -j}^K (P)
$$
This implies the remaining assertions.
\end{proof}

\begin{remark} \label{rem-interp}
Note that the  conjecture of Kalai, Kleinschmidt, and Lee has been shown by giving a
combinatorial interpretation of the first graded Betti numbers of a simplicial polytope.
By duality, it follows that the second last non-trivial graded Betti numbers have a
combinatorial interpretation, too.  However, it is not possible to find combinatorial
interpretations of all graded Betti numbers because, in general, the Betti numbers
depend on the characteristic of the ground field (cf.\ \cite{th-char}, Example 3.3).
\end{remark}
\medskip


\section{Upper bounds for Betti numbers} \label{sec-Betti}

The key to proving the conjecture of Kalai, Kleinschmidt, and Lee has been to identify
the number of missing $i$-simplices as a certain first graded Betti number. The proof
also shows that in order to compute an upper bound for this number in terms of the
$h$-vector of the polytope, we need to know an upper bound for the Betti numbers of
Cohen-Macaulay algebras. The goal of this section is to establish such bounds. Since the
general case does not take more work than the special case of a Cohen-Macaulay algebra,
we will derive upper bounds for the graded Betti numbers of any arbitrary standard graded
$K$-algebra in terms of its Hilbert function.

Throughout this section we denote by $R$ the polynomial ring $K[x_1,\ldots,x_n]$ over an
arbitrary field $K$ with its standard grading where every variable has degree one. $A
\neq 0$ will be a standard graded $K$-algebra $R/I$ where $I \subset R$ is a proper
homogeneous ideal. For a finitely generated graded $R$-module $M = \oplus_{j \in \Z}
[M]_j$, we denote its graded Betti numbers by
$$
\beta^R_{i j} (M) :=  \dim_K [\Tor^R_i (M, K)]_j.
$$
Since the graded Betti numbers of $M$ do not change under field extensions of $K$, we may
and will assume that the field $K$ is infinite.

The Hilbert function of $M$ is the numerical function $h_M: \Z \to \Z, h_M (j) := \dim_K
[M]_j$. The Hilbert functions of  graded $K$-algebras have been completely classified by
Macaulay. In order to state his result we need some notation.

\begin{notation} \label{not-bin-expansion}
(i) We always use the following convention for binomial coefficients: If $a \in \R$ and
$j \in \Z$ then
$$
\binom{a}{j} := \left \{ \begin{array}{ll}
\frac{a (a-1) \cdots (a-j+1)}{j!} & \mif j > 0 \\
1 & \mif j = 0 \\
0 & \mif j < 0.
\end{array} \right.
$$

(ii) Let $b, d$ be positive integers. Then there are uniquely determined integers $m_d >
m_{d-1} > m_s \geq s \geq 1$ such that
$$
b = \binom{m_d}{d} + \binom{m_{d-1}}{d-1} + \ldots + \binom{m_s}{s}.
$$
This is called the {\em $d$-binomial expansion} of $b$. For any integer $j$ we set
$$
b\pp{d, j}  := \binom{m_d + j}{d + j} + \binom{m_{d-1} + j}{d-1 + j } + \ldots +
\binom{m_s + j}{s + j}.
$$

Of particular importance will be the cases where $j = 1$ or $j = -1$. To simplify
notation, we further define
$$
b\pp{d} := b\pp{d, 1}  = \binom{m_d + 1}{d + 1} + \binom{m_{d-1} + 1}{d } + \ldots +
\binom{m_s + 1}{s + 1}
$$
and
$$
b\pl{d} := b\pp{d, -1}  = \binom{m_d - 1}{d - 1} + \binom{m_{d-1}-1}{d - 2} + \ldots +
\binom{m_s - 1}{s - 1}.
$$

(iii) If $b = 0$, then we put $b\pp{d} = b\pl{d} =  b\pp{d, j} := 0$ for all $j, d \in
\Z$.
\end{notation}
\medskip

Recall that a sequence of non-negative integers $\left (h_j\right )_{j \geq 0}$ is called
an {\em O-sequence} if $h_0 = 1$ and $h_{j+1} \leq h_j\pp{j}$ for all $j \geq 1$. Now we
can state Macaulay's characterization of Hilbert functions \cite{Macaulay} (cf.\ also
\cite{stanley}).

\begin{theorem}[Macaulay] \label{thm-Mac}
For a numerical function $h: \Z \to \Z$, the following conditions are equivalent:
\begin{itemize}
\item[(a)] $h$ is the Hilbert function of a standard graded $K$-algebra;
\item[(b)] $h(j) = 0$ if $j < 0$ and $\{h(j)\}_{j \geq 0}$ is an O-sequence.
\end{itemize}
\end{theorem}

For later use we record some formulas for sums involving binomial coefficients.

\begin{lemma} \label{lem-sum-formulas}
For any positive real numbers $a, b$ and every integer $j \geq 0$, there are the
following identities:
\begin{itemize}
\item[(i)] ${\displaystyle \sum_{k = 0}^j (-1)^k \binom{a+k-1}{k} \binom{b}{j-k} =
\binom{b-a}{j} }$;
\item[(ii)] ${\displaystyle \sum_{k = 0}^j \binom{a+k-1}{k} \binom{b+j-k-1}{j-k} =
\binom{a+b+j-1}{j} }$;
\item[(iii)] ${\displaystyle \sum_{k = 0}^j (-1)^k \binom{a+k}{m} \binom{b}{j-k} =
\sum_{k=0}^{m} \binom{a-k-1}{m-k} \binom{b-k-1}{j}}$  \begin{minipage}[t]{3.1cm}
 if \ $0 \leq m
\leq a$ are integers. \end{minipage}
\end{itemize}
\end{lemma}

\begin{proof}
(i) and (ii) are probably standard. In any case, they follow immediately by comparing
coefficients of power series using the identities $(1+x)^{b-a} = (1+x)^{-a} \cdot
(1+x)^b$ and $(1-x)^{-a-b} = (1-x)^{-a} \cdot (1-x)^{-b}$.

To see part (iii), we first use (ii) and finally (i); we get:
\begin{eqnarray*}
\sum_{k = 0}^j (-1)^k \binom{a+k}{m} \binom{b}{j-k} & = & \sum_{k = 0}^j (-1)^k
 \binom{b}{j-k} \cdot \left \{\sum_{i=0}^m \binom{k+i}{i} \binom{a-1-i}{m-i} \right\} \\
& = & \sum_{i=0}^m \binom{a-1-i}{m-i}  \cdot \left \{ \sum_{k = 0}^j (-1)^k
\binom{k+i}{k} \binom{b}{j-k}
\right \} \\
& = & \sum_{k=0}^{m} \binom{a-1-i}{m-i} \binom{b-i-1}{j},
\end{eqnarray*}
as claimed.
\end{proof}

After these preliminaries we are ready to derive bounds for Betti numbers. We begin with
the special case of modules having a $d$-linear resolution. Recall that the graded module
$M$ is said to have a {\em $d$-linear resolution} if it has a graded minimal free
resolution of the form
$$
\ldots \to R^{\beta_i} (-d-i) \to \ldots \to R^{\beta_1} (-d-1) \to R^{\beta_0} (-d) \to
M \to 0.
$$
Here $\beta_i^R (M) = \sum_{j \in \Z} \beta_{i, j}^R (M) := \beta_i$ is the $i$-th total
Betti number of $M$.

\begin{proposition} \label{prop-lin-res}
Let $M \neq 0$ be a graded $R$-module with a $d$-linear resolution. Then, for every $i
\geq 0$, its $i$-th total graded Betti number is
$$
\beta_i^R (M) = \sum_{j=0}^i  (-1)^j \cdot h_M (d+j) \cdot \binom{n}{i-j}.
$$
\end{proposition}

\begin{proof}
We argue by induction on $i$. The claim is clear if $i = 0$. Let $i > 0$. Using the
additivity of vector space dimensions along exact sequences and the induction hypothesis
we get:
\begin{eqnarray*}
\beta_i^R (M) & = & (-1)^i h_M (d+ i) + \sum_{j=0}^{i-1}  (-1)^{i-1 - j} \cdot \beta_j^R
(M)
\binom{n-1+i-j}{i-j} \\
& = & (-1)^i h_M (d+ i) + \\
& & \sum_{j=0}^{i-1} (-1)^{i-1 - j} \cdot \binom{n-1+i-j}{i-j} \cdot \left \{
\sum_{k=0}^j (-1)^j \cdot  h_M (d+k) \binom{n}{j-k} \right \}\\
& = & (-1)^i h_M (d+ i) + \\
&& \sum_{k=0}^{i-1} (-1)^{k} \cdot h_M (d+k) \cdot \left \{ \sum_{j=k}^{i-1} (-1)^{i-1-j}
 \binom{n-1-i-j}{i-j} \binom{n}{j-k} \right \}\\
& = & (-1)^i h_M (d+ i) + \\ && \sum_{k=0}^{i-1} (-1)^{k} \cdot h_M (d+k) \cdot \left \{
\sum_{j=1}^{i-k} (-1)^{j-1}
 \binom{n+j-1}{j} \binom{n}{i-k-j} \right \}\\
& = & (-1)^i h_M (d+ i) + \sum_{k=0}^{i-1} (-1)^{k} \cdot h_M (d+k) \cdot \binom{n}{i-k}
\end{eqnarray*}
according to Lemma \ref{lem-sum-formulas}(i). Now the claim follows.
\end{proof}

It is amusing and useful to apply this result to a case where we know the graded Betti
numbers.

\begin{example} \label{ex-powers}
Consider the ideal $I = (x_1,\ldots,x_n)^d$ where $d>0$. Its minimal free resolution is
given by an Eagon-Northcott complex. It has a $d$-linear resolution and its Betti numbers
are (cf., e.g., the proof of \cite{MN3}, Corollary 8.14):
$$
\beta_i^R (I) = \binom{d+i-1}{i} \binom{n+d-1}{d+i}.
$$
Since the Hilbert function of $I$ is, for all $j \geq 0$, $h_I (d+j) =
\binom{n+d+j-1}{d+j}$, a comparison with Proposition \ref{prop-lin-res} yields:
\begin{equation} \label{eq-sum}
\binom{d+i-1}{i} \binom{n+d-1}{d+i}  =  \sum_{j=0}^i (-1)^j \cdot \binom{n+d+j-1}{d+j}
\binom{n}{i-j}.
\end{equation}
\end{example}
\medskip

Now we will compute the graded Betti numbers of lex-segment ideals. Recall that an ideal
$I \subset R$ is called a {\em lex-segment} ideal if, for every $d$, the ideal $I\dd{d}$
is generated by the first $\dim_k [I]_d$ monomials in the lexicographic order of the
monomials in $R$. Here $I\dd{d}$ is the ideal that is generated by all the polynomials of
degree $d$ in $I$. For every graded $K$-algebra $A= R/I$ there is a unique lex-segment
ideal $I^{lex} \subset R$ such that $A$ and $R/I^{lex}$ have the same Hilbert function.
For further information on lex-segment ideals we refer to \cite{BH-book}.

\begin{lemma} \label{lem-lin-lex-seg}
Let $I \subset R$ be a proper lex-segment ideal whose generators all have degree $d$.
Consider the $d$-binomial expansion of $b := h_{R/I} (d)$:
$$
b = \binom{m_d}{d} + \binom{m_{d-1}}{d-1} + \ldots + \binom{m_s}{s}.
$$
Then the Betti numbers of $A := R/I$ are for all $i \geq 0$:
\begin{eqnarray*}
\beta^R_{i+1} (A) & = & \beta^R_{i+1, i+d} (A)\\
& = & \binom{n+d-1}{d+i} \binom{d+i-1}{d-1} - \sum_{k=s}^d \sum_{j=0}^{m_k - k}
\binom{m_k -j -1}{k-1} \binom{n-1-j}{i}.
\end{eqnarray*}
(Note that  according to Notation \ref{not-bin-expansion}, the sum on the right-hand side
is zero if $b = 0$.)
\end{lemma}

\begin{proof}
Gotzmann's Persistence Theorem \cite{gotzmann} implies that the Hilbert function of $A$
is, for $j \geq 0$, $h_A (d+j) = b\pp{d, j}$ and that $I$ has a $d$-linear resolution.
Hence Proposition \ref{prop-lin-res}  in conjunction with Formula (\ref{eq-sum}) and
Lemma \ref{lem-sum-formulas}(iii) provides:
\begin{eqnarray*}
\beta_{i+1}^R (A) & = & \beta_i^R (I) = \sum_{j=0}^i (-1)^j \cdot h_I(d+j) \cdot
\binom{n}{i-j} \\
& = & \sum_{j=0}^i (-1)^j \left [\binom{n+d+j-1}{d+j} - b\pp{d, j} \right ] \cdot
\binom{n}{i-j} \\
& = & \binom{n+d-1}{d+i} \binom{d+i-1}{i} - \sum_{j=0}^i (-1)^j \cdot \left [\sum_{k=s}^d
\binom{m_k + j}{k+j} \right ] \cdot \binom{n}{i-j} \\
& = &  \binom{n+d-1}{d+i} \binom{d+i-1}{i} -  \sum_{k=s}^d \left [\sum_{j=0}^i (-1)^j
\cdot
\binom{m_k + j}{m_k - k}  \cdot \binom{n}{i-j} \right ]\\
& = &  \binom{n+d-1}{d+i} \binom{d+i-1}{i} -  \sum_{k=s}^d \sum_{j=0}^{m_k - k}
\binom{m_k - j -1}{k-1}  \cdot \binom{n-1-j}{i},
\end{eqnarray*}
as claimed.
\end{proof}

The above formulas simplify in the extremal cases.

\begin{corollary} \label{cor-linear-lex-seg}
Adopt the notation and assumptions of Lemma \ref{lem-lin-lex-seg}. Then
\begin{itemize}
\item[(a)] ${\displaystyle \beta^R_1 (A) = \binom{n+d-1}{d} - b}$;
\item[(b)] ${\displaystyle \beta^R_n (A) = \binom{n+d-2}{d-1} - b\pl{d} }$.
\end{itemize}
\end{corollary}

\begin{proof}
Part (a) being clear, we restrict ourselves to show (b). Since $\binom{n-1-j}{n-1} = 0$
for $j > 0$, Lemma \ref{lem-lin-lex-seg} immediately gives
$$
\beta^R_n (A) = \binom{n+d-2}{d-1} - \sum_{k=s}^d \binom{m_k -1}{k-1} =
\binom{n+d-2}{d-1} -  b\pl{d}.
$$
\end{proof}

Now, we can compute the non-trivial graded Betti numbers of an arbitrary lex-segment
ideal.

\begin{proposition} \label{prop-res-lex-segment}
Let $I \subset R$ be an arbitrary proper lex-segment ideal and let $d \geq 2$ be an
integer. Set $A := R/I$ and consider the $d$-binomial expansion
$$
h_A (d) =: \binom{m_d}{d} + \binom{m_{d-1}}{d-1} + \ldots + \binom{m_s}{s}
$$
and the $(d-1)$-binomial expansion
$$
h_A (d-1) =: \binom{n_{d-1}}{d-1} + \binom{n_{d-2}}{d-2} + \ldots + \binom{n_t}{t}.
$$
Then we have for all $i \geq 0$:
$$
\beta^R_{i+1, i+d} (A) = \beta_{i+1, i+d} (h_A, n)
$$
where
$$
\beta_{i+1, i+d} (h_A, n) := \sum_{k=t}^{d-1} \sum_{j=0}^{n_k - k} \binom{n_k - j}{k}
\binom{n-1-j}{i} - \sum_{k=s}^{d} \sum_{j=0}^{m_k - k} \binom{m_k - 1 - j}{k-1}
\binom{n-1-j}{i}.
$$
\end{proposition}

\begin{proof}
As noticed above, since $I$ is a lex-segment ideal, for every $j \in \Z$, the ideal
$I\dd{j}$ has a $j$-linear resolution, i.e.\ the ideal $I$ is componentwise linear. Hence
\cite{HH}, Proposition 1.3, gives for all $i \geq 0$:
\begin{equation} \label{eq-HH-formula}
\beta^R_{i+1, i+d} (A) = \beta^R_{i+1} (R/I\dd{d}) - \beta^R_{i+1} (R/\fm I\dd{d-1})
\end{equation}
where $\fm = (x_1,\ldots,x_n)$ is the homogeneous maximal ideal of $R$.

Since $I\dd{d-1}$ is generated in degree $d-1$, the ideals $I\dd{d-1}$ and $\fm
I\dd{d-1}$ have the same Hilbert function in all degrees $j \geq d$. Thus, using the
assumption $d \geq 2$, Gotzmann's Persistence Theorem (\cite{gotzmann}) provides:
$$
h_{R/\fm I\dd{d-1}} (d-1+j) = h_{R/I\dd{d-1}} (d-1+j) = h_{A} (d-1)\pp{d-1, j} \quad
\fall j \geq 1.
$$
It is easy to see that  $\fm I\dd{d-1}$ has a $d$-linear resolution because $I\dd{d-1}$
has a $(d-1)$-linear resolution. Hence, as in the proof of Lemma \ref{lem-lin-lex-seg},
Proposition \ref{prop-lin-res} provides
$$
\beta^R_{i+1} (R/\fm I\dd{d-1})  =  \binom{n+d-1}{d+i} \binom{d+i-1}{d-1} -
\sum_{k=t}^{d-1} \sum_{j=0}^{n_k - k} \binom{n_k -j}{k} \binom{n-1-j}{i}.
$$
Plugging this and the result of Lemma \ref{lem-lin-lex-seg} into the Formula
(\ref{eq-HH-formula}), we get our claim.
\end{proof}

Again, the formula simplifies in the extremal cases. We will use the result in the
following section.

\begin{corollary} \label{cor-lex-seg}
Adopt the notation and assumptions of Proposition \ref{prop-res-lex-segment}. Then:
\begin{itemize}
\item[(a)] ${\displaystyle \beta^R_{1, d} (A) = \beta_{1, d} (h_A, n) = h_A (d-1)\pp{d-1}
- h_A (d)}$;
\item[(b)] ${\displaystyle \beta^R_{n, n-1+d} (A) = \beta_{n, n-1+d} (h_A, n) = h_A (d-1)
- (h_A (d))\pl{d}}$.
\end{itemize}
\end{corollary}

\begin{proof}
This follows from the formula given in Proposition \ref{prop-res-lex-segment}.
\end{proof}

In Proposition \ref{prop-res-lex-segment} we left out the  case  $d \leq 1$ which is easy
to deal with. We need:

\begin{definition} \label{def-betas}
Let $h$ be the Hilbert function of graded $K$-algebra such that $h(1) \leq n$. Then we
define, for all integers $i \geq 0$ and $d$, the numbers $ \beta_{i+1, i+d} (h, n)$ as in
Proposition \ref{prop-res-lex-segment} if $d \geq 2$ and otherwise:
$$
\beta_{i+1, i+d} (h, n) := \left \{ \begin{array}{cl}
\binom{n - h(1)}{i+1} & \mif d = 1 \\[1ex]
0 & \mif d \leq 0.
\end{array} \right.
$$
Moreover, if $i \leq 0$ we set:
$$
\beta_{i, j} (h, n) := \left \{ \begin{array}{cl}
1 & \mif (i, j) = (0, 0) \\[1ex]
0 & {\rm otherwise}.
\end{array} \right.
$$
\end{definition}

\begin{lemma} \label{lem-small-d}
Let $A = R/I \neq 0$ be any graded $K$-algebra. Then we have for all integers $i, d$ with
 $d \leq 1$:
$$
\beta^R_{i+1, i+d} (A) =  \beta_{i+1, i+d} (h_A, n).
$$
\end{lemma}

\begin{proof}
Since $A$ has as $R$-module just one generator in degree zero, this is clear if $d \leq
0$. Furthermore, $I\dd{1}$ is generated by a regular sequence of length $n - h_A (1)$.
Its minimal free resolution is given by the Koszul complex. Hence, the claim follows for
$d = 1$ because $\beta^R_{i+1, i+1} (A) = \beta^R_{i, i+1} (I\dd{1})$.
\end{proof}

Combined with results of Bigatti, Hullet, and Pardue, we get the main result of this
section: bounds for the graded Betti numbers of a $K$-algebra as an $R$-module in terms
of its Hilbert function and the dimension of $R$.

\begin{theorem} \label{thm-Betti-bounds}
Let $A = R/I \neq 0$ be a graded $K$-algebra. Then its graded Betti numbers are bounded
by
$$
\beta^R_{i+1, i+j} (A) \leq  \beta_{i+1, i+j} (h_A, n) \quad (i, j \in \Z).
$$

Furthermore, equality is attained for all integers $i, j$ if $I$ is a lex-segment ideal.
\end{theorem}

\begin{proof}
Let $I^{lex} \subset R$ be the lex-segment ideal such that $A$ and $R/I^{lex}$ have the
same Hilbert function. Then we have for all integers $i, j$ that
$$
\beta^R_{i+1, i+j} (A) \leq  \beta^R_{i+1, i+j} (R/I^{lex})
$$
according to Bigatti \cite{bigatti} and Hulett \cite{hulett} if $\chara K = 0$  and to
Pardue \cite{pardue} if $K$ has positive characteristic. Since Proposition
\ref{prop-res-lex-segment} and  Lemma \ref{lem-small-d} yield
$$
\beta^R_{i+1, i+j} (R/I^{lex}) = \beta_{i+1, i+j} (h_A, n) \quad (i, j \in \Z),
$$
our claims follow.
\end{proof}

\begin{remark} \label{rem-Hilb-syz}
Note that Theorem \ref{thm-Betti-bounds} gives in particular that $\beta^R_{i+1,i+ d} (A)
= 0$  if $i \geq n$, in accordance with Hilbert's Syzygy Theorem.
\end{remark}
\medskip

We conclude this section by discussing the graded Betti numbers of Cohen-Macaulay
algebras with the so-called \WLP.

Let $A = R/I$ be a graded Cohen-Macaulay $K$-algebra of Krull dimension $d$ and let
$l_1,\ldots,l_d \in [R]_1$ be sufficiently general linear forms. Then $\overline{A} :=
A/(l_1,\ldots,l_d) A$ is called the {\em Artinian reduction} of $A$. Its Hilbert function
and graded Betti numbers as module over $\overline{R} := R/(l_1,\ldots,l_d) R$ do not
depend on the choice of the forms $l_1,\ldots,l_d$. The Hilbert function of
$\overline{A}$ takes positive values in only finitely many degrees. The sequence of these
positive integers $\underline{h} = (h_0, h_1,\ldots, h_r)$ is called the {\em $h$-vector}
of $A$. We set $\beta_{i+1, i+d} (\underline{h}, n-d) := \beta_{i+1, i+d}
(h_{\overline{A}}, n-d)$. Using this notation we get.

\begin{corollary} \label{cor-CM-alg}
Let $A = R/I$ be a Cohen-Macaulay graded $K$-algebra of dimension $d$ with $h$-vector
$\underline{h}$. Then its graded Betti numbers satisfy
$$
\beta^R_{i+1, i+j} (A) \leq  \beta_{i+1, i+j} (\underline{h}, n-d) \quad (i, j \in \Z).
$$
\end{corollary}

\begin{proof}
If $l \in [R]_1$ is a not a zerodivisor of $A$, then the graded Betti numbers of $A$ as
$R$-module agree with the graded Betti numbers of $A/l A$ as $R/l R$ module (cf., e.g.,
\cite{MN3}, Corollary 8.5). Hence, by passing to the Artinian reduction of $A$, Theorem
\ref{thm-Betti-bounds} provides the claim.
\end{proof}

\begin{remark} \label{rem-vanishing} Note that, for any O-sequence $\underline{h} =
(1, h_1,\ldots, h_r)$ with $h_r > 0$, Definition \ref{def-betas} provides $\beta_{i+1,
i+j} (\underline{h}, m) = 0$ for all $i, m \geq 0$ if $j \leq 1$ or $j \geq r+2$.
\end{remark}

Recall that an Artinian graded $K$-algebra $A$  has the so-called \WLP\ if there is an
element $l \in A$ of degree one such that, for each $j \in \Z$, the multiplication
$\times l: [A]_{j-1} \to [A]_j$ has maximal rank. The Cohen-Macaulay $K$-algebra $A$ is
said to have the {\em \WLP} if its Artinian reduction has the \WLP.

\begin{remark} \label{rem-WLP}
The Hilbert functions of Cohen-Macaulay algebras with the \WLP\ have been completely
classified in \cite{wlp}, Proposition 3.5. Moreover, Theorem 3.20 in \cite{wlp} gives
optimal upper bounds on their graded Betti numbers in terms of the Betti numbers of
certain lex-segment ideals. Thus, combining this result with Theorem
\ref{thm-Betti-bounds}, one gets upper bounds for the Betti numbers of these algebras in
terms of their Hilbert functions. In general, these bounds are strictly smaller than the
bounds of Corollary \ref{cor-CM-alg} for Cohen-Macaulay algebras that do not necessarily
have the \WLP.
\end{remark}
\medskip

The $h$-vectors of graded Gorenstein algebras with the \WLP\ are precisely the
SI-sequences (cf.\ \cite{MN3}, Theorem 6.3, or \cite{harima}, Theorem 1.2). For their
Betti numbers we obtain:

\begin{theorem} \label{thm-gor-betti}
Let $\underline{h} = (1,h_1,\dots,h_u,\dots, h_r)$ be an SI-sequence where $h_{u-1 } <
h_u = \cdots = h_{r-u} > h_{r-u+1}$.  Put $\underline{g} =
(1,h_1-1,h_2-h_1,\dots,h_u-h_{u-1})$.  If $A = R/I$ is a Gorenstein graded $K$-algebra of
dimension $d$ with the \WLP\ and $h$-vector $\underline{h}$, then its graded Betti
numbers satisfy
$$
\beta^R_{i+1, i+j} (A) \leq \left \{
\begin{array}{ll}
\beta_{i+1, i+j} (\underline{g}, m) & \hbox{if $j \leq r-u$} \\
\beta_{i+1, i+j} (\underline{g}, m) + \beta_{g_1-i,r+h_1-i-j} (\underline{g}, m) &
\hbox{if $r-u+1 \leq j \leq u+1$} \\
\beta_{g_1-i,r+h_1-i-j} (\underline{g}, m) & \hbox{if $j \geq u+2$}
\end{array}
\right.
$$
where $m := n-d-1 = \dim R - d- 1$.
\end{theorem}

\begin{proof}
This follows immediately by combining \cite{MN3}, Theorem 8.13, and Theorem
\ref{thm-Betti-bounds}.
\end{proof}


\section{Explicit bounds for the number of missing simplices} \label{sec-missing}

We now return to the consideration of simplicial polytopes. To this end we will
specialize the results of  Section \ref{sec-Betti} and then discuss some applications.

We begin by simplifying somewhat our notation. Let $P$ be a simplicial $d$-polytope with
$f$-vector $\underline{f}$. It is well-known that the $h$-vector of the Stanley-Reisner
ring $K[P]$ agrees with the $h$-vector of $P$ as defined in the introduction.
Furthermore, in Section \ref{sec-conj} we defined the graded Betti numbers of $K[P]=
R/I_P$ by resolving $K[P]$ as an $R$-module where $R$ is a polynomial ring of dimension
$f_0$ over $K$, i.e.\
$$
\beta_{i, j}^K (P) = \beta_{i, j}^R (K[P]).
$$
Note that the Stanley-Reisner ideal $I_P$ does not contain any linear forms. The graded
Betti numbers of $P$ agree with the graded Betti of the Artinian reduction of $K[P]$ as
module over a polynomial ring of dimension $f_0 - d = h_1$. Thus, we can simplify the
statements of the bounds of $\beta_{i, j}^K (P)$ by setting:

\begin{notation} \label{not-betas}
Using the notation introduced above Corollary \ref{cor-CM-alg} we define for every
O-sequence $\underline{h}$
$$
\beta_{i+1, i+j} (\underline{h}) := \beta_{i+1, i+j} (\underline{h}, h_1).
$$
\end{notation}
\medskip

Notice that $\beta_{i+1, i+j} (\underline{h}) = 0$ if $i\geq 0$ and $j \leq 1$.

In this section we will primarily use the $g$-vector of a polytope which is defined as
follows:

\begin{definition} \label{deg-g-vector}
Let $P$ be a simplicial polytope with $h$-vector $\underline{h} := (h_0,\ldots, h_d)$.
Then the g-Theorem (Theorem \ref{g-thm}) shows that there is a unique integer $u$ such
that  $h_{u-1 } < h_u = \cdots = h_{d-u} > h_{d-u+1}$.  The vector $\underline{g} =
(g_0,\ldots, g_u) := (1,h_1-1,h_2-h_1,\dots,h_u-h_{u-1})$ is called the {\em $g$-vector}
of $P$. All its entries are positive.
\end{definition}

Some observations are in order.

\begin{remark} \label{rem-g-vec}
(i) By its definition, the $g$-vector of the polytope $P$ is uniquely determined by the
$h$-vector of $P$. The g-Theorem shows that the $h$-vector of $P$ (thus also its
$f$-vector) can be recovered from its $g$-vector, provided the dimension of $P$ is given.

(ii) The g-Theorem also gives an estimate of the length of the $g$-vector because it
implies $2 u \leq d = \dim P$.
\end{remark}

Now we can state our explicit bounds for the Betti numbers of a polytope.

\begin{theorem} \label{thm-poly-betti}
Let $K$ be a field of characteristic zero and let  $\underline{g} = (g_0,\dots,g_u)$ be
an O-sequence with $g_u > 0$. Then we have:
\begin{itemize}
\item[(a)] If $P$ is a simplicial $d$-polytope with $g$-vector $\underline{g}$, then:
$$
\beta^K_{i+1, i+j} (P) \leq \left \{
\begin{array}{ll}
\beta_{i+1, i+j} (\underline{g}) & \hbox{if $j \leq d-u$} \\
\beta_{i+1, i+j} (\underline{g}) + \beta_{g_1-i,d+h_1-i-j} (\underline{g}) &
\hbox{if $d-u+1 \leq j \leq u+1$} \\
\beta_{g_1-i,d+g_1+1-i-j} (\underline{g}) & \hbox{if $j \geq u+2$.}
\end{array}
\right.
$$
\item[(b)] In {\rm (a)} equality is attained  for all integers $i, j$ if $P$ is the
$d$-dimensional Billera-Lee polytope with $g$-vector $\underline{g}$.
\end{itemize}
\end{theorem}

\begin{proof}
According to Stanley \cite{stanley2} (cf.\ also \cite{McM-simple}), the Stanley-Reisner
ring of every simplicial polytope has the \WLP.  Hence part (a) is a consequence of
Theorem \ref{thm-gor-betti}. Part (b) follows from \cite{MN3}, Theorem 9.6, and Theorem
\ref{thm-Betti-bounds}, as pointed out in the proof of Theorem \ref{thm-betti}.
\end{proof}

We have seen in Section \ref{sec-conj} that the number of empty $j$-simplices of the
simplicial polytope $P$ is equal to the  Betti number $\beta^K_{1, j+1} (P)$. Thus, we
want to make the preceding bounds more explicit if $i = 0$. At first, we treat a trivial
case.

\begin{remark} \label{rem-trivial-g1}
Notice that the $g$-vector has length one, i.e.\ $u = 0$ if and only if the polytope $P$
is a simplex. In this case, its Stanley-Reisner ideal is a principal ideal generated by a
monomial of degree $d = \dim P$.
\end{remark}
\medskip

In the following result we stress when the Betti numbers vanish. Because of Remark
\ref{rem-trivial-g1}, it is harmless to assume that $u \geq 1$. We use Notation
\ref{not-bin-expansion}.

\begin{corollary} \label{cor-generators}
Let  $\underline{g} = (g_0,\dots,g_u)$ be an O-sequence with $g_u > 0$ and $u \geq 1$.
Set $g_{u+1} := 0$. Then we have:
\begin{itemize}
\item[(a)] If $P$ is a simplicial $d$-polytope with $g$-vector $\underline{g}$, then
there are the following bounds:
\begin{itemize}
\item[(i)] If $d \geq 2 u + 1$, then
$$
\beta^K_{1, j} (P) \leq \left \{
\begin{array}{ll}
g_{j-1}\pp{j-1} - g_j & \mif \ 2 \leq j \leq u+1 \\[.5ex]
g_{d+1-j} - (g_{d+2-j})\pl{d+2-j} & \mif \ d-u+1 \leq j \leq d \\[.5ex]
0 & {otherwise; }
\end{array}
\right.
$$
\item[(ii)] If $d = 2 u$, then
$$
\beta^K_{1, j} (P) \leq \left \{
\begin{array}{ll}
g_{j-1}\pp{j-1} - g_j & \mif \ 2 \leq j \leq u \\[.5ex]
g_u\pp{u} + g_u & \mif \ j = u+1 \\[.5ex]
g_{d+1-j} - (g_{d+2-j})\pl{d+2-j} & \mif \ u+2 \leq j \leq d \\[.5ex]
0 & {otherwise; }
\end{array}
\right.
$$
\end{itemize}
\item[(b)] In {\rm (a)} equality is attained  for all integers $j$ if $P$ is the
$d$-dimensional Billera-Lee polytope with $g$-vector $\underline{g}$.
\end{itemize}
\end{corollary}

\begin{proof}
Since the first Betti numbers of any polytope do not depend on the characteristic of the
field, the claims follow from Theorem \ref{thm-poly-betti} by taking into account
Corollary \ref{first-betti}, Corollary \ref{cor-lex-seg}, and the fact that $\beta_{i+1,
i+j}(\underline{g}) = 0$ if $i \geq 0$ and either $j \leq 1$ or $j \geq u+2$ due to
Remark \ref{rem-vanishing}.
\end{proof}

To illustrate the last result, let us consider an easy case.

\begin{example} \label{ex-ci}
Let $P$ be a simplicial $d$-polytope with $g_1 = 1$. Then its Stanley-Reisner ideal $I_P$
is a Gorenstein ideal of height two, thus a complete intersection. Indeed, since the
$g$-vector of $P$ is an O-sequence, it must be $\underline{g} = (g_0,\ldots, g_u) =
(1,\ldots,1)$. Hence Corollary \ref{cor-generators} provides that $I_P$ has exactly two
minimal generators, one of degree $u+1$ and one of degree $d-u+1$. Equivalently, $P$ has
exactly two empty simplices, one of dimension $u$ and one of dimension $d-u$.
\end{example}
\medskip

As an immediate consequence of Corollary \ref{cor-generators} we partially recover
\cite{Kalai-94}, Proposition 3.6.

\begin{corollary} \label{cor-no-empty}
Every simplicial $d$-polytope has no empty faces of dimension $j$ if $u+1 \leq j \leq
d-u-1$.
\end{corollary}

\begin{remark} \label{rem-converse}
Kalai's Conjecture 8 in \cite{Kalai-94} states that the following converse of Corollary
\ref{cor-no-empty} should be true: If there is an integer $k$ such that $d \geq 2k$ and
the simplicial $d$-polytope has no empty simplices of dimension $j$ whenever $k \leq j
\leq d-k$, then $u < k$. Kalai has proved this if  $k=2$ in \cite{Kalai-87}. Our results
provide the following weaker version of Kalai's conjecture:

If there is an integer $k$ such that $d \geq 2k$ and {\em every} simplicial $d$-polytope
with $g$-vector $(g_0,\ldots,g_u)$ has no empty simplices of dimension $j$ whenever $k
\leq j \leq d-k$, then $u < k$.

Indeed, this follows by the sharpness of the bounds in Corollary \ref{cor-generators}.
\end{remark}

Now we want to make some existence results of Kalai and Perles effective.  As
preparation, we state:

\begin{corollary} \label{cor-up-to-k}
Let $P$ be a simplicial $d$-polytope with $g$-vector $\underline{g} = (g_0,\dots,g_u)$
where $u \geq 1$. Set $g_{u+1} = 0$. Then the number $N(k)$ of empty simplices of $P$
whose dimension is at most $k$, is bounded above as follows:
\begin{equation*}
N(k) \leq \left \{ \begin{array}{ll}

g_1 + \sum_{j=1}^k \left \{ g_j\pp{j} - g_j \right \} - g_{k+1} & \mif \ 1 \leq k \leq
\min
\{u, d-u-1 \}; \\[1ex]
N(u) & \mif \ u < k < d-u \\[1ex]
\begin{minipage}{7.1cm}
$g_1 + g_{d-k}\pp{d-k} + \sum_{j=1}^{d-k-1} \left \{ g_j\pp{j} - g_j \right \}  +
\\
\hspace*{2.1cm} \sum_{j = d-k+1}^u \left \{ g_j\pp{j} - (g_j)\pl{j} \right \}$
 \end{minipage}
 & \mif \ d-u
\leq k < d
\end{array} \right.
\end{equation*}

Furthermore, for each $k$, the bound is attained if $P$ is the Billera-Lee $d$-polytope
with $g$-vector $\underline{g}$.
\end{corollary}

\begin{proof}
By Corollary \ref{cor-no-empty}, this is clear if $u < k < d-u$. In any case, we know
that $N(k) = \sum_{j=2}^{k+1} \beta^K_{1, j} (P)$. Thus, using Corollary
\ref{cor-generators} carefully, elementary calculations provide the claim. We omit the
details.
\end{proof}

The last result immediately gives:

\begin{corollary} \label{cor-total}
If $P$ is a simplicial polytope with $g$-vector $\underline{g} = (g_0,\dots,g_u)$ where
$u \geq 1$, then its total number of empty simplices is at most
$$
\binom{g_1+ 2}{2} - 1 + \sum_{j = 2}^u \left \{ g_j\pp{j} - (g_j)\pl{j} \right \}.
$$

Furthermore, this bound is attained if $P$ is any Billera-Lee polytope with  $g$-vector
$\underline{g}$.
\end{corollary}

\begin{proof}
Use Corollary \ref{cor-up-to-k} with $k = d-1$ and recall that $g_1\pp{1} =
\binom{g_1+1}{2}$.
\end{proof}

\begin{remark} It is somewhat surprising that the bound in Corollary \ref{cor-total} does {\em not} depend
on the dimension of the polytope. In contrast, the other bounds (cf., e.g., Corollary
\ref{cor-up-to-k}) do depend on the dimension $d$ of the polytope.
\end{remark}

In view of Corollary \ref{cor-up-to-k}, the following elementary facts will be useful.

\begin{lemma} \label{lem-monotonie}
Let  $k$ be a positive integer. If $a \geq b$ are non-negative integers, then
\begin{itemize}
\item[(a)] \hspace*{3cm} ${\displaystyle a\pp{k} - a\pl{k} \geq  b\pp{k} - b\pl{k}; }$\\[-1ex]
\item[(b)] \hspace*{3.3cm} ${\displaystyle a\pp{k} - a \geq  b\pp{k} - b; }$\\[-1ex]
\item[(c)] \hspace*{4.1cm} ${\displaystyle  a\pl{k} \geq   b\pl{k}. }$
\end{itemize}
\end{lemma}

\begin{proof}
We show only (a). The proofs of the other claims are similar and only easier.

To see (a), we begin by noting, for integers $m \geq j
>0$, the identity
\begin{equation} \label{eq-comp}
\binom{m+1}{j+1} - \binom{m-1}{j-1} = \binom{m}{j+1} + \binom{m-1}{j}.
\end{equation}
Now we use induction on $k \geq 1$. Since $a\pp{1} - a\pl{1} = \binom{a+1}{2} - 1$, the
claim is clear if $k = 1$. Let $k \geq 2$. Consider the $k$-binomial expansions
$$
a =: \binom{m_k}{k} + \binom{m_{k-1}}{k-1} + \ldots + \binom{m_s}{s} \quad {\rm and}
\quad b =: \binom{n_{k}}{k} + \binom{n_{k-1}}{k-1} + \ldots + \binom{n_t}{t}.
$$
Since $a \geq b$, we get $m_k \geq n_k$. We distinguish two cases.
\smallskip

{\em Case 1}: Let $m_k = n_k$. Then the claim follows by applying the induction
hypothesis to $a - \binom{m_k}{k} \geq b - \binom{m_k}{k}$.
\smallskip

{\em Case 2}: Let $m_k > n_k$. Using $n_i \leq n_k-k+i$ and Formula (\ref{eq-comp}), we
get
\begin{eqnarray*}
b\pp{k} - b\pl{k} &= & \sum_{i = t}^{k} \left \{ \binom{n_i}{i+1} + \binom{n_i - 1}{i}
\right \} \\
& \leq & \sum_{i = 1}^{k} \left \{ \binom{n_k - k + i}{i+1} + \binom{n_k - k - 1 + i}{i}
\right \} \\
& = & \binom{n_k+1}{k+1} + \binom {n_k}{k} - (n_k - k + 2) \\
& < & \binom{m_k}{k+1} + \binom{m_k - 1}{k}
\end{eqnarray*}
because $n_k < m_k$. The claim follows since Formula (\ref{eq-comp}) gives
$\binom{m_k}{k+1} + \binom{m_k - 1}{k} \leq a\pp{k} - a\pl{k}$.
\end{proof}

\begin{remark}
In general, it is not true that $a > b$ implies $a\pp{k} - a\pl{k} > b\pp{k} - b\pl{k}$.
For example, if $k \geq 2$ and $a-1 = b= \binom{m}{k} > 0$, then $a\pp{k} - a\pl{k} =
b\pp{k} - b\pl{k}$.
\end{remark}
\medskip

We are ready to establish optimal bounds that depend only on the dimension and the number
of vertices.

\begin{theorem} \label{thm-kal-2-7} Let $P$ be a simplicial $d$-polytope with $d + g_1 +
1$ vertices which is not a simplex. Then there is the following bound on the number
$N(k)$ of empty simplices of $P$ whose dimension is $\leq k$:
\begin{equation*}
N(k) \leq \left \{ \begin{array}{ll}
\binom{g_1 + k}{g_1 - 1}  & \mif \ 1 \leq k < \frac{d}{2}; \\[1ex]
\binom{g_1 + \left \lfloor \frac{d}{2} \right \rfloor}{g_1 - 1} + \binom{g_1 + \left
\lfloor \frac{d}{2} \right \rfloor - 1}{g_1 - 1} & \mif \ \frac{d}{2} \leq k< d.
\end{array} \right.
\end{equation*}

Furthermore, for each $k$, the bound is attained if $P$ is the Billera-Lee $d$-polytope
with $g$-vector $(g_0,\ldots,g_u)$ where $g_j = \binom{g_1 + j - 1}{j}$, $0 \leq j \leq
u$, and $u = \min \{k, \left \lfloor \frac{d}{2} \right \rfloor \}$.
\end{theorem}

\begin{proof}
Let  $\underline{g} = (g_0,\dots,g_u)$ be the $g$-vector of $P$. Since $P$ is not a
simplex, we have $u \geq 1$. We have to distinguish two cases.
\smallskip

{\em Case 1}: Let $k <  \frac{d}{2}$. If $k > u$, then we formally set $g_{u+1} = \ldots
= g_{\left \lfloor \frac{d}{2} \right \rfloor} =0$. Since $k < \frac{d}{2} \leq d - u$,
Corollary \ref{cor-up-to-k} provides:
$$
N(k) \leq g_1 + \sum_{j=1}^k \left \{ g_j\pp{j} - g_j \right \} - g_{k+1}.
$$
According to Lemma \ref{lem-monotonie}, the sum on the right-hand side becomes maximal if
$g_2,\ldots,g_k$ are as large as possible and $g_{k+1} = 0$. The latter means $u = k$.
Macaulay's Theorem \ref{thm-Mac} implies $g_j \leq \binom{g_1 + j - 1}{j}$. Now an easy
computation provides the bound in this case. It is sharp because $(g_0,\ldots,g_k)$,
where $g_j = \binom{g_1 + j - 1}{j}$, is a $g$-vector of a simplicial $d$-polytope by the
$g$-Theorem, thus Corollary \ref{cor-up-to-k} applies.
\smallskip

{\em Case 2}: Let $\frac{d}{2} \leq k < d$. First, let us also assume that $k \geq d-u$.
Then Corollary \ref{cor-up-to-k} gives:
$$
N(k) \leq g_1 + g_{d-k}\pp{d-k} + \sum_{j=1}^{d-k-1} \left \{ g_j\pp{j} - g_j \right \} +
\sum_{j = d-k+1}^u \left \{ g_j\pp{j} - (g_j)\pl{j} \right \}.
$$
Again, Lemma \ref{lem-monotonie} shows that, for fixed $u$, the bound is maximized if
$g_j = \binom{g_1 + j - 1}{j}$, $0 \leq j \leq u$. This provides
$$
N(k) \leq \binom{g_1 + u}{g_1 -1} +\binom{g_1 + u - 1}{g_1 - 1}.
$$
Since $u \leq \frac{d}{2}$, our bound follows in this case.

Second, assume $k < d-u$. Then $u \leq \frac{d}{2} \leq k < d-u$ yields $u <
\frac{d}{2}$. Thus Corollary \ref{cor-up-to-k} provides $N(k) = N(u)$, but $N(u) \leq
\binom{g_1 +u}{g_1 - 1}$ by Case 1. This concludes the proof of the bound in Case 2. Its
sharpness is shown as in Case 1.
\end{proof}

As immediate consequence we obtain:

\begin{corollary} \label{cor-kalai-2-7}
Every simplicial polytope, which is not a simplex, has at most $\binom{g_1 + k}{g_1 - 1}
+ \binom{g_1 + k-1}{g_1 - 1}$ empty simplices of dimension $\leq k$.
\end{corollary}

\begin{remark} \label{rem-comp-bounds}
Kalai \cite{Kalai-94}, Theorem 2.7, has first given an estimate as in Corollary
\ref{cor-kalai-2-7}. His bound is
$$
N(k) \leq (g_1 + 1)^{k+1} \cdot (k+1)!.
$$
Comparing with our bound, we see that Kalai's bound is asymptotically not optimal for
$g_1 \gg 0$.
\end{remark}
\medskip

Notice that the bound on $N(k)$ in Theorem \ref{thm-kal-2-7} does not depend on $k$ if $k
\geq \frac{d}{2}$. This becomes plausible by considering cyclic polytopes.

\begin{example}
(i)
 Recall that a {\em cyclic polytope} $C(f_0, d)$ is a $d$-dimensional simplicial
polytope which is the convex hull of $f_0$ distinct points on the moment curve
$$
\{ (t, t^2,\ldots,t^d) \s t \in \mathbb{R} \}.
$$
Its combinatorial type depends only on $f_0$ and $d$.

According to McMullen's Upper Bound Theorem (\cite{McM-cyclic}),  the cyclic polytope
$C(f_0, d)$ has  the maximal $f$-vector among all simplicial $d$-polytopes with $f_0$
vertices. Theorem \ref{thm-kal-2-7} shows that it also has the maximal total number of
empty simplices among these polytopes. Indeed, this follows by comparing with the main
result in \cite{th-cyclic} (cf.\ also \cite{MN3}, Corollary 9.10) which provides that
$C(f_0, d)$ has $\binom{g_1 + \left \lfloor \frac{d}{2} \right \rfloor}{g_1 - 1} +
\binom{g_1 + \left \lfloor \frac{d}{2} \right \rfloor - 1}{g_1 - 1}$ empty simplices.
Moreover, the empty simplices of $C(f_0, d)$ have either dimension $\frac{d}{2}$ if $d$
is even or dimensions $\frac{d-1}{2}$ and $\frac{d+1}{2}$ if $d$ is odd. This explains
why the bound on $N(k)$ in Theorem \ref{thm-kal-2-7} does not change if $k \geq
\frac{d}{2}$.

(ii) If $P$ is a simplicial $d$-polytope with $f_0 \geq d + 2$ vertices, then Theorem
\ref{thm-kal-2-7} gives for its number of empty edges
$$
N(1) \leq \left \{ \begin{array}{ll}
\frac{f_0 (f_0 - 3)}{2} & \mif \ d = 2 \\[1ex]
\binom{f_0 - d}{2} & \mif \ d \geq 3.
\end{array} \right.
$$
If $d=2$, the bound is always attained because $\frac{f_0 (f_0 - 3)}{2}$ is the number of
``missing diagonals'' of a convex $f_0$-gon.
\end{example}
\medskip

\begin{remark} \label{rem-skeleta}
Recall that the {\em $k$-skeleton} of an arbitrary $d$-polytope $P$ is the set of all
faces of $P$ whose dimension is at most $k$. Perles \cite{perles} has shown:

The number of combinatorial types of $k$-skeleta of $d$-polytopes with $d+g_1+1$ vertices
is bounded by a function in $k$ and $g_1$.

Following \cite{Kalai-94}, the proof of this result can be reduced to the case where the
polytopes are simplicial. Then one concludes by using  a bound on $N(k)$ because the
$k$-skeleton of a simplicial polytope is determined by its set of empty simplices of
dimension $\leq k$.
\end{remark}
\medskip

In \cite{Kalai-94} Kalai sketches an argument showing that the number of empty simplices
can be bounded with very little information on the $g$-vector. Below, we will slightly
correct \cite{Kalai-94}, Theorem 3.8, and give explicit bounds. We use Notation
\ref{not-bin-expansion}.

\begin{theorem} \label{thm-kalai-3-8}
Fix  integers $j \geq k \geq 1$ and $b \geq 0$. Let $P$ be a simplicial polytope
$d$-polytope $P$ with $g_k \leq b$ where we define $g_i = 0$ if $i > u$. If $d \geq j +
k$, then the number of empty $j$-simplices  of $P$ is bounded by
$$
\left \{
\begin{array}{ll}
b\pp{k, j-k+1} & \mif \ j < \frac{d}{2} \\
b\pp{k, j-k+1} + b\pp{k, j-k} & \mif \ j = \frac{d}{2} \\
b\pp{k, d-j-k} & \mif \ j > \frac{d}{2} \\
\end{array} \right.
$$
\end{theorem}

\begin{proof}
We have to bound $\beta_{1,j+1}^K (K[P])$.  By Corollary \ref{cor-no-empty}, $P$ has no
empty $j$-simplices if $u+1 \leq j \leq d-u-1$. Thus, we may assume that $1 \leq j \leq
u$ or $d-u \leq j \leq d-1$.
\smallskip

{\em Case 1}: Assume $1 \leq j \leq u \leq \frac{d}{2}$. Then Corollary
\ref{cor-generators} provides if $j <   \frac{d}{2}$:
$$
\beta_{1,j+1}^K (K[P]) \leq g_j\pp{j} - g_{j+1}.
$$
Using Lemma \ref{lem-monotonie}, we see that the bound is maximized if $g_{j+1} = 0$ and
$g_j$ is as large as possible. Since the $g$-vector is an O-sequence, we get $g_j \leq
g_k\pp{k, j-k} \leq b\pp{k, j-k}$. Our claimed bound follows.

If $j = \frac{d}{2}$, then we get $j = u = \frac{d}{2}$. Hence Corollary
\ref{cor-generators} gives:
$$
\beta_{1,j+1}^K (K[P]) \leq g_j\pp{j}  + g_j.
$$
Now the bound is shown as above.
\smallskip

{\em Case 2}: Assume $\frac{d}{2} \leq d-u \leq j \leq d-1$. By the above considerations,
we may also assume that $j \neq \frac{d}{2}$. Thus, Corollary \ref{cor-generators}
provides:
$$
\beta_{1,j+1}^K (K[P]) \leq g_{d-j} - (g_{d+1-j})\pl{d+1-j}.
$$
Using our assumption $d-j \geq k$, we conclude as above.
\end{proof}

\begin{remark}
(i) In  \cite{Kalai-94}, Theorem 3.8, the existence of bounds as in the above result is
claimed without assuming $d \geq j+k$. However, this is impossible, as Case 2 in the
above proof shows. Indeed, if $d-j < k$ and $d > j > \frac{d}{2}$, then knowledge of
$g_k$ does not give any information on $g_{d-j}$. In particular,  $g_{d-j}$ can be
arbitrarily large preventing the existence of a bound on $\beta_{1,j+1}^K (K[P])$ in
terms of $g_k, j, k$ in this case.

For a somewhat specific example, fix $k = j = 2$ and $d=3$. Then the Billera-Lee
3-polytope with $g$-vector $(1, g_1)$ has $g_1$ empty 2-simplices.

(ii) Note that the bounds in Theorem \ref{thm-kalai-3-8} are sharp if $g_k = b$. This
follows from the proof.
\end{remark}

If we only know that $d$ is large enough compared to $j$ and $k$, then we have the
following weaker bound.

\begin{corollary} \label{cor-j-simp}
Fix  integers $j \geq k \geq 1$,  $b \geq 0$, and $d \geq j+k$. Then the number of empty
$j$-simplices  of every simplicial $d$-polytope with $g_k \leq b$ is at most $b\pp{k,
j-k+1} + b\pp{k, j-k}$.
\end{corollary}

\begin{proof}
By Theorem \ref{thm-kalai-3-8}, it remains to consider the case where $j > \frac{d}{2}$.
But then $d-j < j$, thus $b\pp{k, d-j-k} \leq b\pp{k, j-k}$, and we conclude again by
using Theorem \ref{thm-kalai-3-8}.
\end{proof}

\begin{remark}
Notice that the bound in Corollary \ref{cor-j-simp} is independent of the number of
vertices of the polytope and its dimension, provided the latter is large enough.
\end{remark}

In essence, all the bounds on the number of empty simplices are bounds on certain first
graded Betti numbers of the Stanley-Reisner ring of a simplicial polytope. As such, using
Theorem \ref{thm-gor-betti}, they can be extended to bounds for the first graded Betti
numbers of any graded Gorenstein algebra with the Weak Lefschetz property. We leave this
and analogous considerations for higher Betti numbers to the interested reader.
\bigskip

\noindent
{\bf Acknowledgments}
\smallskip

The author would like to thank Gil Kalai, Carl Lee, and Juan Migliore for
  motivating discussions, encouragement, and  helpful comments.


\end{document}